\documentclass[leqno,draft,11pt]{article}
\usepackage{amssymb,amsmath,latexsym,theorem,a4wide}

\def\sameenum{}

\ifx\mydefs\undefined\else \fi
\let\mydefs\relax

\allowdisplaybreaks[2]

\ifx\loadcyr\undefined\else
\input cyracc.def
\makeatletter
 at 1\@ptsize pt
 at 1\@ptsize pt
 at 1\@ptsize pt
\makeatother

\fi

\def\gobble#1{}
\def\fixsup#1#2{{#1\let\dp\gobble\mathstrut}^#2_}
\let\mr\mathrel

\def\bme{\hskip.75em\relax}

\let\eq\leftrightarrow
\let\EQ\Leftrightarrow

\def\iff{\quad\text{iff}\quad}
\let\ET\bigwedge
\let\TO\Rightarrow
\def\?{\mathbin?}

\let\model\vDash
\let\nmodel\nvDash

\newbox\circlebox
\setbox\circlebox\hbox{$\bigcirc$}
\def\circled#1{%
  \setbox0\hbox to\wd\circlebox{\hss$#1$\hss}\wd0=0pt
  \box0\copy\circlebox}

\let\fii\varphi
\let\tet\vartheta
\let\ep\varepsilon

\def\greek#1{$\expandafter\greeknum\csname c@#1\endcsname$}
\makeatletter
\def\greeknum#1{\ifcase#1\or\alpha\or\beta\or\gamma\or\delta\or\ep
      \or\digamma\or\zeta\or\eta\or\tet\or\iota\else\@ctrerr\fi}
\makeatother

\def\p#1{\langle#1\rangle}

\def\lh#1{\lvert#1\rvert}

\let\sset\subseteq
\let\nsset\nsubseteq

\let\Sset\supseteq

\let\onto\twoheadrightarrow

\def\pw#1{\mathcal P(#1)}

\newcommand\rpair[3][3em]{\mathrel{%
   \begin{matrix}%
     \strut\smash{\xrightonto{\hbox to#1{\hss$#2$\hss}}}\\[-1.7ex]%
     \strut\smash{\xleftembed[\hbox to#1{\hss$#3$\hss}]{}}%
   \end{matrix}}}
\makeatletter
\newcommand\xrightonto[2][]{\ext@arrow 0359\rightontofill{#1}{#2}}
\newcommand\xleftembed[2][]{\ext@arrow 3095\leftembedfill{#1}{#2}}
\def\leftembedfill{\arrowfill@\leftarrow\relbar\hookleftnoarrow}
\def\rightontofill{\arrowfill@\relbar\relbar\onto}
\def\hookleftnoarrow{\DOTSB\relbar\joinrel\rhook}
\makeatother

\mathchardef\#="2023 

\ifx\busspf\undefined\else
\usepackage{bussproofs}
\EnableBpAbbreviations

\fi

\let\dia\Diamond
\def\diadot{\centerdot\dia}
\ifx\symlasy\undefined
  
  \def\boxdot{{\origboxdot}}
  \def\centerdot#1{{%
    \setbox0\hbox{$\mathop{#1}$}\dimen0 \ht0
    \setbox0\hbox{$#1$}\advance\dimen0 -\ht0
    \setbox2\hbox to\wd0{\hss$\mathop{\cdot}$\hss}\wd2=0pt
    \lower\dimen0\box2\box0 }}
\else
  \def\boxdot{\centerdot\Box}
  \def\centerdot#1{{%
     \setbox0\hbox{$#1$}%
     \raise0.206\ht0\hbox to\wd0{\hss$\cdot$\hss}%
     \kern-\wd0 \box0 }}
\fi

\let\sls|

\def\Up{{\setbox0\hbox{$\uparrow$}%
         \lower\dp0\hbox to\wd0{\hss\vrule width4pt height.4pt\hss}%
         \kern-\wd0\box0}}
\def\UP{{\setbox0\hbox{$\uparrow$}%
         \lower\dp0\hbox to\wd0{\hss\vrule width4pt height.4pt\hss}%
         \kern-\wd0\copy0\kern-\wd0\raise.35ex\box0}}
\def\Down{{\setbox0\hbox{$\downarrow$}%
         \raise\ht0\hbox to\wd0{\hss\vrule width4pt depth.4pt\hss}%
         \kern-\wd0\box0}}

\newif\ifnadm

\def\doadm{\mathrel{%
   \setbox0 \hbox{$\mathop\vdash$}\dimen0 \ht0
   \setbox0 \hbox{$\vdash$}\advance\dimen0 -\ht0
   \vrule width.8\fontdimen8 \textfont3 height\ht0 depth\dp0
   \mkern-1mu
   \lower\dimen0 \hbox{$\vcenter{%
      \ifnadm
        \setbox0 \hbox{$\scriptstyle\sim\mathstrut$}%
        \hbox{\hbox to\wd0{\hss$\scriptstyle/$\hss}\kern-\wd0 \box0 }%
      \else
        \hbox{$\scriptstyle\sim\mathstrut$}%
      \fi}$}}}

\def\nrstyle#1#2#3{%
  \setbox0\hbox{$#1\bigcirc$}%
  \vcenter{\hbox to\wd0{\hss$#2#3$\hss}}%
  \kern-\wd0\box0 }

\DeclareMathOperator\id{id}

\DeclareMathOperator\Sub{Sub}

\def\st{\expandafter\hat}

\let\lgc\mathbf

\def\I{{\bullet}}
\def\R{{\circ}}

\mathcode`\*="0203

\mathchardef\mhyphen="2D

\def\noproof{\leavevmode\unskip\bme\vadjust{}\nobreak\hfill$\qed$\par}
\let\qed\Box
\newenvironment{Pf}[1][]
  {\par\noindent\textit{Proof\optpar{#1}:}\bme\ignorespaces}
  {\noproof\pagebreak[2]\vskip\medskipamount\ignorespacesafterend}
\def\optpar#1{\ifx\relax#1\relax\else\ #1\fi}
\def\qedhere{\relax\ifmmode\eqno\qed\expandafter\aftergroup
                   \else\noproof\fi\noqed}
\def\noqed{\let\noproof\relax}

\theoremstyle{plain}
\ifx\shortthm\undefined
\newtheorem{Thm}{Theorem}[section]
\else
\newtheorem{Thm}{Theorem}
\fi
\newtheorem{Prop}[Thm]{Proposition}
\newtheorem{Cor}[Thm]{Corollary}

\newtheorem{Fact}[Thm]{Fact}
\newtheorem{Obs}[Thm]{Observation}

\newtheorem{Ass}[Thm]{Assumption}
\newtheorem{Cl}{Claim}[Thm]
\ifx\shortthm\undefined
\def\theCl{\arabic{Cl}}
\fi

\theorembodyfont\upshape
\newtheorem{Def}[Thm]{Definition}

\newtheorem{Exm}[Thm]{Example}

\newenvironment{Pf*}{\let\qed\qedCl\Pf}\endPf

\usepackage[reftex]{theoremref}
\newif\iflinenumbers
\linenumberstrue

{\catcode`\^^I=13 \catcode`\^^M=13
\gdef\doalgo#1#2\end#{\hbox to\hsize{\hss \let^^I\qquad%
  \def\\^^M{\nobreak\hfil\break\vadjust{}\qquad}%
  \fboxsep1em \linenum0 %
  \fbox{\hsize#1\vbox{%
  \everypar{\advance\linenum1 %
      \hbox to1.2em{%
           \hss\iflinenumbers$\scriptstyle\the\linenum$\hskip.6em\fi}}%
  #2}}\hss}\end}}
\newcount\linenum

\def\key{\relax\ifmmode\expandafter\mathbf\else\expandafter\textbf\fi}

\def\allowhyphens{\nobreak\hskip0pt\relax}

\DeclareRobustCommand*\magiclparen{\ifmmode(\else\textup(\allowhyphens\fi}
\DeclareRobustCommand*\magicrparen{\ifmmode)\else\textup)\fi}
\let\lparen=(  \let\rparen=)
\def\magicparon{\catcode`\(\active\catcode`\)\active}
\def\magicparoff{\catcode`\(12 \catcode`\)12 }
\AtBeginDocument{\ifx\ifPreview\iftrue\else\magicparon\fi}
\magicparon
\let (=\magiclparen  \let )=\magicrparen

\ifx\sameenum\undefined
  
  \ifx\enumup\undefined
    
  \else
    
  \fi

\else

\fi

\magicparoff

\mathchardef\comma=\mathcode`\,
{\catcode`\,=\active \gdef,{\comma\penalty\relpenalty}}

\providecommand\dedic{\thanks{Supported by
grant IAA100190902 of GA AV \v CR, Center of Excellence CE-ITI under the grant
P202/12/G061 of GA \v CR, and RVO: 67985840.}}
\author{Emil Je\v r\'abek\dedic\\[\medskipamount]
Institute of Mathematics of the Academy of Sciences\\
\small \v Zitn\'a 25,
115\:67 Praha 1,
Czech Republic,
email: \texttt{jerabek@math.cas.cz}
}

\def\bdi#1{#1^{\boxdot^{-1}}}
\def\wbc{BDP}
\def\sbc{SBDP}
\let\frm\mathcal

\title{Cluster expansion and the boxdot conjecture}

\begin{document}
\maketitle

\begin{abstract}
The boxdot conjecture asserts that every normal modal logic that
faithfully interprets~$\lgc T$ by the well-known boxdot translation is
in fact included in~$\lgc T$. We confirm that the conjecture is
true. More generally, we present a simple semantic condition on modal
logics~$L_0$ which ensures that the largest logic where $L_0$
embeds faithfully by the boxdot translation is $L_0$ itself. In particular, this
natural generalization of the boxdot conjecture holds for $\lgc{S4}$,
$\lgc{S5}$, and $\lgc{KTB}$ in place of~$\lgc T$.
\end{abstract}

\section{The boxdot translation}\label{sec:boxdot}
The \emph{boxdot translation} is the mapping $\fii\mapsto\fii^\boxdot$ from the
language of monomodal logic into itself that preserves
propositional variables, commutes with Boolean connectives, and
satisfies
\[(\Box\fii)^\boxdot=\boxdot\fii^\boxdot,\]
where $\boxdot\fii$ is an abbreviation for~$\fii\land\Box\fii$. It is
easy to see that for any normal modal logic~$L$, the set of formulas
interpreted in~$L$ by the boxdot translation,
\[\bdi L=\{\fii:{}\vdash_L\fii^\boxdot\},\]
is also a normal modal logic (nml), and contains the logic $\lgc
T=\lgc K\oplus\Box p\to p$. The boxdot translation is a faithful
interpretation of~$\lgc T$ in the smallest nml~$\lgc K$ (i.e.,
$\bdi{\lgc K}=\lgc T$), and more generally, in any logic between $\lgc
K$ and~$\lgc T$. The \emph{boxdot conjecture}, formulated by French and
Humberstone~\cite{fr-hum}, states that the converse also holds:
\[\bdi L=\lgc T\implies L\sset\lgc T.\]
French and Humberstone proved the conjecture for logics~$L$
axiomatized by formulas of modal degree~$1$, and
Steinsvold~\cite{steinsv} has shown it for logics of the form $L=\lgc
K\oplus\dia^h\Box^ip\to\Box^j\dia^kp$, but the full conjecture
remained unsettled.

In this paper we will establish the boxdot conjecture. The argument
actually uses only one particular property of~$\lgc T$, namely that
Kripke frames for~$\lgc T$ can be blown up by duplicating each node in
a frame in a certain way, and this allows us to state an analogue of
the conjecture for a large class of reflexive logics.

Given an arbitrary nml~$L_0\Sset\lgc T$, we may ask about the
structure of logics~$L$ in which $L_0$ faithfully embeds by the boxdot
translation ($\bdi L=L_0$). First, there is always at least one such logic, for
example $L_0$ itself: this follows from the observation that
$\vdash_\lgc T\fii\eq\fii^\boxdot$ for every formula~$\fii$. If
$L_0=\lgc K\oplus X$ for some set of axioms~$X$, then the logic
$L_0^\boxdot=\lgc K\oplus\{\fii^\boxdot:\fii\in X\}$ has the property
\[L_0\sset\bdi L\iff L_0^\boxdot\sset L,\]
in particular $L_0^\boxdot\sset L_0$ is the \emph{smallest} logic in which
$L_0$ faithfully embeds. Clearly, the set of logics such that $\bdi
L=L_0$ is \emph{convex}: if $L_1\sset L_2\sset L_3$ and
$\bdi{L_1}=\bdi{L_3}=L_0$, then $\bdi{L_2}=L_0$. Finally, if
$\{L_c:c\in C\}$ is a chain of logics linearly ordered by inclusion
such that $\bdi{L_c}=L_0$ for each~$c\in C$, the logic
$L=\bigcup_{c\in C}L_c$ also satisfies $\bdi L=L_0$. It follows from
Zorn's lemma that every logic in which $L_0$ faithfully embeds is
included in a \emph{maximal} such logic.

Thus, if $\{L_m:m\in M\}$ is the set of all maximal logics such that
$\bdi{L_m}=L_0$, the set of all logics in which $L_0$ faithfully
embeds by the boxdot translation consists of the union
\begin{equation}\label{eq:1}
\bigcup_{m\in M}[L_0^\boxdot,L_m]
\end{equation}
of intervals in the lattice of normal modal logics. Notice that $L_0$
is itself one of the maximal logics~$L_m$, as
$\bdi L=L$ for any $L\Sset L_0$ (or $L\Sset\lgc T$ for that matter).
For a nontrivial example, $\lgc{A^*}=\lgc{GL}\oplus\Box\Box
p\to\Box(\boxdot p\to q)\lor\Box(\boxdot q\to p)$ is a maximal logic
in which~$\lgc{S4Grz}$ embeds~\cite[Exer.~9.26]{cha-zax}.

The original boxdot conjecture states that $\lgc T$ is the largest
logic~$L$ such that $\bdi L=\lgc T$. In accordance with this, we define
that a nml $L_0\Sset\lgc T$ has the \emph{boxdot
property}, if
\[\tag{\wbc} \bdi L=L_0\implies L\sset L_0\]
holds for every nml~$L$. In light of the discussion above, \wbc\
for~$L_0$ is equivalent to the claim that there is only one maximal
logic~$L$ such that $\bdi L=L_0$, or in other words, that the union
in~\eqref{eq:1} reduces to the single interval $[L_0^\boxdot,L_0]$.

What we are going to show is that \wbc\ holds for all logics~$L_0$
satisfying a natural semantic condition (which we will define
precisely in Section~\ref{sec:main-result}). Since the condition applies to~$\lgc
T$, this also establishes the original boxdot conjecture. Moreover,
our proof shows that under the same condition, $L_0$ has the
\emph{strong boxdot property}:
\[\tag{\sbc} \bdi L\sset L_0\implies L\sset L_0.\]
We do not know whether \wbc\ and~\sbc\ are equivalent in general, though
they of course express the same condition when~$L_0=\lgc T$.

\section{Preliminaries}\label{sec:preliminaries}
We first recall elementary concepts from relational semantics to agree on the
notation.
\begin{Def}\th\label{def:kripke}
A \emph{Kripke frame} is a pair $\frm W=\p{W,R}$, where $R$ is a binary
relation on a set~$W$. A \emph{model} based on the frame~$\frm W$ is a
triple $\frm M=\p{W,R,{\model}}$, where the \emph{valuation~$\model$} is
a relation between elements of~$W$ and modal formulas, written as
$\frm M,w\model\fii$, which satisfies
\begin{align*}
\frm M,w&\model\fii\to\psi\iff
    \frm M,w\nmodel\fii\text{ or }\frm M,w\model\psi,\\
\frm M,w&\model\Box\fii\iff
    \forall v\in W\,(w\mr Rv\TO\frm M,v\model\fii),
\end{align*}
and similarly for other Boolean connectives.
A formula~$\fii$ \emph{holds} in~$\frm M$ if $\frm M,w\model\fii$ for every~$w\in W$,
and it is \emph{valid} in~$\frm W$, written as $\frm W\model\fii$, if
$\frm M\model\fii$ for every model~$\frm M$ based on~$\frm W$. If $L$
is a nml, $\frm W$ is a \emph{Kripke $L$-frame}
if $\frm W\model\fii$ for every~$\fii\in L$. A logic $L$ is
\emph{sound} wrt a class~$C$ of frames if every $\frm W\in C$ is an
$L$-frame, and it is \emph{complete} wrt~$C$ if $\fii\notin L$ implies
$\frm W\nmodel\fii$ for some~$\frm W\in C$.

A \emph{cluster} in a transitive frame~$\p{W,R}$ is an equivalence
class of the equivalence relation~$\sim$ on~$W$ defined by
\[w\sim v\iff w=v\lor(w\mr Rv\land v\mr Rw).\]
\end{Def}

We will also consider the general frame semantics. The reason is
mostly esthetic---we do not want to encumber our results with the
arbitrary restriction to Kripke-complete logics~$L_0$ which does not seem to
have anything to do with the problem under investigation. A reader who
is only interested in the original boxdot conjecture ($L_0=\lgc T$),
or more generally in \wbc\ or~\sbc\ for Kripke-complete logics~$L_0$, may
safely ignore general frames in what follows.
\begin{Def}\th\label{def:genfr}
A \emph{general frame} is a triple $\frm W=\p{W,R,A}$, where $\p{W,R}$
is a Kripke frame, and $A$ is a family of subsets of~$W$ which is
closed under Boolean operations, and under the operation
\[\Box X=\{w\in W:\forall v\in W\,(w\mr Rv\TO v\in X)\}.\]
Sets~$X\in A$ are called \emph{admissible}. A model $\frm
M=\p{W,R,{\model}}$ is based on~$\frm W$ if the set
\[\{w\in W:\frm M,w\model p\}\]
is admissible for every variable~$p$ (which implies the same holds for all
formulas). A formula is valid in~$\frm W$ if it holds in all models
based on~$\frm W$, and the notions of $L$-frames, soundness, and
completeness are defined accordingly. A Kripke frame~$\p{W,R}$ can be
identified with the general frame~$\p{W,R,\pw W}$.
\end{Def}

We will also need basic validity-preserving operations on frames and models.
\begin{Def}\th\label{def:frameops}
A Kripke frame~$\p{W',R'}$ is a \emph{generated subframe} of~$\frm W=\p{W,R}$
if $W'\sset W$, and $R'=R\cap(W'\times W)$. (Note that $W'\sset W$ is
a carrier of a generated subframe of~$\frm W$ iff it is upward closed
under~$R$.) If $\p{W,R}$ and $\p{W',R'}$ are Kripke frames, a mapping
$f\colon W\to W'$ is a \emph{p-morphism}, provided
\begin{enumerate}
\item $w\mr Rv$ implies $f(w)\mr{R'}f(v)$,
\item if $f(w)\mr{R'}v'$, there is $u\in W$ such that $w\mr Ru$ and~$f(u)=v'$,
\end{enumerate}
for every $w,v\in W$ and~$v'\in W'$. Notice that the image of a
p-morphism is always a generated subframe of~$\frm{W'}$.

Similarly, a general frame~$\p{W',R',A'}$ is a generated subframe of a
frame~$\p{W,R,A}$ if the Kripke frame~$\p{W',R'}$ is a generated
subframe of~$\p{W,R}$, and $A'=\{X\cap W':X\in A\}$. A p-morphism from
a general frame~$\p{W,R,A}$ to~$\p{W',R',A'}$ is a p-morphism from the
Kripke frame $\p{W,R}$ to~$\p{W',R'}$ which additionally satisfies
\begin{enumerate}
\setcounter{enumi}2
\item $f^{-1}[X']\in A$ for every~$X'\in A'$.
\end{enumerate}
\end{Def}
\begin{Fact}[{{\cite[Thm.~3.14, Prop.~5.72]{brv}}}]\th\label{fact:pres}
Let $\frm W$ and~$\frm{W'}$ be Kripke or general frames, and $\fii$ a
formula valid in~$\frm W$. If $\frm{W'}$ is a generated subframe
of~$\frm W$, or if there exists a p-morphism from~$\frm W$
onto~$\frm{W'}$, then $\frm{W'}\model\fii$.
\noproof\end{Fact}
In more detail, p-morphisms preserve truth in models in the following
way.
\begin{Fact}[{{\cite[Prop.~2.14]{brv}}}]\th\label{fact:modpres}
If $\frm M=\p{W,R,{\model}}$ and~$\frm M'=\p{W',R',{\model}}$ are
models, $\fii$ is a formula, and $f$ is a p-morphism from~$\p{W,R}$ to~$\p{W',R'}$ such
that $\frm M,w\model p_i\EQ\frm M',f(w)\model p_i$ for every $w\in
W$ and every variable~$p_i$ occurring in~$\fii$, then $\frm
M,w\model\fii\EQ\frm M',f(w)\model\fii$.
\noproof\end{Fact}

We have been using $\vdash_L\fii$ as a synonym for~$\fii\in L$,
however we will also extend this notation to allow for nonlogical
axioms.
\begin{Def}\th\label{def:vdash}
For any nml~$L$, $\vdash_L$ denotes the global consequence
relation of~$L$: if $X$ is a set of formulas, and $\fii$ a formula,
then $X\vdash_L\fii$ iff $\fii$ has a finite derivation using elements
of~$X$, theorems of~$L$, modus ponens, and necessitation.
\end{Def}
The global consequence relation satisfies the following version of the
deduction theorem, where
$\Box^n\fii=\underbrace{\Box\cdots\Box}_{n\text{ boxes}}\fii$, and
$\Box^{\le n}\fii=\ET_{i=0}^n\Box^i\fii$.
\begin{Fact}[{{\cite[Thm.~3.51]{cha-zax}}}]\th\label{fact:dt}
$X\vdash_L\fii$ iff there is a finite subset $X_0\sset X$
and a natural number~$n$ such that $\Box^{\le n}\ET X_0\to\fii\in L$.
\noproof\end{Fact}

\section{Motivation}\label{sec:motivation}
Before we introduce the semantic condition that will guarantee the
\sbc\ for a logic~$L_0\Sset\lgc T$, let us give some intuition.
This section is mostly informal, its purpose is to explain that the
condition does not
come out of blue, but follows naturally from the properties of the
boxdot translation. However, the argument uses a somewhat heavier
machinery than the rest of the paper, hence it is intended for
readers familiar with the structure theory of transitive modal logics
(see e.g.~\cite{cha-zax} for more background).
It can be skipped without losing continuity, though some of the
counterexamples in \th\ref{exm:bdpfail} may be worth bearing in mind.

Assume that $L\nsset L_0$, we would like to show that $\bdi
L\nsset L_0$. If $L_0$ has the finite model property, there is a
finite rooted $L_0$-frame~$\frm F$ which is not an $L$-frame. If $L$
and~$L_0$ are transitive (i.e., extensions of~$\lgc{K4}$), we can
consider a frame formula~$A_\frm F$ as in Fine~\cite{fine:ff}:
$A_\frm F$ is invalid in a transitive frame~$\frm W$ iff $\frm
F$ is a p-morphic image of a generated subframe of~$\frm W$. (Note
that we will not actually use frame formulas in the proof of our main
result below.) Clearly $\nvdash_{L_0}A_\frm F$.

Do we have $\vdash_LA_\frm F^\boxdot$? Well, if not, and if $L$
is Kripke-complete for simplicity, there is an $L$-frame $\frm
W=\p{W,R}$  such that $\frm W\nmodel A_\frm F^\boxdot$. This
means that $A_\frm F$ is invalid in the reflexivization $\frm
W^\R=\p{W,R\cup\id}$, hence $\frm F$ is a p-morphic image of a
generated subframe of~$\frm W^\R$. Since the class of $L$-frames is
closed under generated subframes, we may assume that there is a
p-morphism~$f$ from $\frm W^\R$~itself onto~$\frm F$.

If $f$ were a p-morphism from~$\frm W$ to~$\frm F$, then $\frm F$
would be an $L$-frame contrary to the assumptions. This contradiction
would show that $\vdash_LA_\frm F^\boxdot$, providing an
example for~$\bdi L\nsset L_0$. In general, $f$ does not have to be a p-morphism
from~$\frm W$ to~$\frm F$, however the only
way this can fail is that for some~$w\in W$, there is no~$u\in W$ such
that $w\mr Ru$ and~$f(w)=f(u)$. The key observation is that this
cannot happen if all clusters of~$\frm F$ have more than one element:
then we can fix $v'\ne f(w)$ in the same cluster as~$f(w)$; since $f$
is a p-morphism from~$\frm W^\R$ to~$\frm F$, there must be $w\mr
Rv\mr Ru$ such that $f(v)=v'$ and~$f(u)=f(w)$ (where necessarily $w\ne
v\ne u$), and we have~$w\mr Ru$ by transitivity.

This suggests that a (transitive) logic~$L_0$ will satisfy \sbc\ if the
class of $L_0$-frames is closed under the operation of blowing up each
cluster by adding new elements. On the other hand, there are examples
implying that some condition of that sort is necessary, which shows
that we are on the right track:
\begin{Exm}\th\label{exm:bdpfail}
The logic~$\lgc{S4Grz}$ corresponds to noetherian partially ordered
frames. In particular, such frames only have one-element clusters.
$\lgc{S4Grz}$ does not have the \wbc: it is well known that
$\bdi{\lgc{GL}}=\lgc{S4Grz}$, but $\lgc{GL}\nsset\lgc{S4Grz}$ (in
fact, $\lgc{GL}$ and~$\lgc{S4Grz}$ have no consistent common
extension).

Similarly, top-most clusters in $\lgc{S4.1}$-frames can only have one
element, and \wbc\ duly fails for~$\lgc{S4.1}$: for example,
$\bdi L=\lgc{S4.1}$, where $L=\lgc{K4}\oplus\diadot\Box\bot$
(this logic is complete wrt finite 
transitive frames whose top clusters are irreflexive).

The motivational argument above may give the false impression that
it is enough for (S)\wbc\ if we can blow up one-element clusters in
$L_0$-frames to have at least two elements, and leave the rest unchanged,
but this is only an artifact of the stipulation that $L$ is
transitive. For a simple counterexample, let $L_0$ be the logic
of the two-element cluster~$\frm C_2$, and $L$ the logic of the
(nontransitive) frame $\frm I_2=\p{\{0,1\},\{\p{0,1},\p{1,0}\}}$. Then
$\bdi L=L_0$ as $\frm I_2^\R=\frm C_2$, but $L\nsset L_0$: e.g.,
$\vdash_Lp\land\Box(\Box p\to p)\to\Box p$, which is not valid
in~$L_0$.
\end{Exm}

\section{The boxdot property}\label{sec:main-result}
The concept of clusters does not make much sense for nontransitive
frames, nevertheless the operation of doubling the size of each
cluster can be easily generalized to arbitrary frames while retaining
its most salient properties.
\begin{Def}\th\label{def:dbl}
For any Kripke frame $\frm W=\p{W,R}$, we define a new frame $2\frm W=\p{2W,2R}$ by putting $2W=W\times\{0,1\}$, and
\[\p{w,a}\mr{2R}\p{v,b}\iff w\mr Rv\]
for any $w,v\in W$ and $a,b\in\{0,1\}$.
If $\frm W=\p{W,R,A}$ is a general frame, we put $2\frm
W=\p{2W,2R,2A}$, where
\[2A=\bigl\{(X\times\{0\})\cup(Y\times\{1\}):X,Y\in A\bigr\}.\]
\end{Def} 
We remark that a similar construction of modal frames from
intuitionistic frames is employed in~\cite{cha-zax}
under the name~$\tau_2$ for investigation of the G\"odel translation.

The following is immediate from the definition.
\begin{Obs}\th\label{obs:2w}
If $\frm W$ is a Kripke or general frame, the natural projection
$\pi\colon2W\to W$ is a p-morphism.
\noproof\end{Obs}
Motivated by the discussion in Section~\ref{sec:motivation}, we will
consider logics with the following property.
\begin{Ass}\th\label{ass:dbl}
$L_0$ is a normal modal logic complete with respect to a class~$C$ of Kripke
or general frames such that $2\frm W$ is an $L_0$-frame for
every~$\frm W\in C$.
\end{Ass}
Notice that even though it is not explicitly demanded, the assumption
implies that $L_0$ is
sound wrt~$C$ because of \th\ref{fact:pres,obs:2w}.
\begin{Exm}\th\label{exm:ass}
It is readily seen that if $\frm W$ is reflexive, symmetric, or
transitive, then so is~$2\frm W$. Thus, \th\ref{ass:dbl} holds for
$\lgc T$, $\lgc{KTB}$, $\lgc{S4}$, and $\lgc{S5}$. For reflexive
transitive frames, $2\frm W$ has the simple geometric interpretation
of expanding each cluster of~$\frm W$ to twice its original size, which implies that
\th\ref{ass:dbl} also holds for $\lgc{S4.2}$ and~$\lgc{S4.3}$.
\end{Exm}

We come to our main result.
\begin{Thm}\th\label{thm:main}
Every logic $L_0\Sset\lgc T$ satisfying \th\ref{ass:dbl} has the
strong boxdot property.
\end{Thm}
\begin{Pf}
Let $L$ be a nml such that $L\nsset L_0$. We fix a formula~$\fii\in L$
such that $\fii\notin L_0$, and we will construct a formula~$\chi$
such that $\chi^\boxdot\in L$ and~$\chi\notin L_0$, witnessing that
$\bdi L\nsset L_0$.

Let $\Sub(\fii)$ denote the set of all subformulas of~$\fii$. We
choose a propositional variable $q\notin\Sub(\fii)$, and consider
the finite set of formulas
\[X=\bigl\{\Box(q^e\to\psi)\to\psi:
   e\in\{0,1\},\Box\psi\in\Sub(\fii)\bigr\},\]
where $q^1=q$, $q^0=\neg q$.
\begin{Cl}\th\label{cl:1}\
\begin{enumerate}
\item\label{item:4}
$X^\boxdot\vdash_\lgc K\Box\psi^\boxdot\to\psi^\boxdot$ for
every $\Box\psi\in\Sub(\fii)$ (even without the use of necessitation).
\item\label{item:5}
$X^\boxdot\vdash_\lgc K\psi\eq\psi^\boxdot$ for every $\psi\in\Sub(\fii)$.
\end{enumerate}
\end{Cl}
\begin{Pf*}
\ref{item:4}: We have
\[\vdash_\lgc K 
  q^e\land\Box\psi^\boxdot\to\boxdot(q^{1-e}\to\psi^\boxdot)\]
for $e=0,1$, and $X^\boxdot$ includes the formula
$\boxdot(q^{1-e}\to\psi^\boxdot)\to\psi^\boxdot$, which yields
\[X^\boxdot\vdash_\lgc K q^e\to(\Box\psi^\boxdot\to\psi^\boxdot).\]
The result follows using $\vdash_\lgc K q\lor\neg q$.

\ref{item:5}: By induction on the complexity of~$\psi$. The steps for
variables and Boolean connectives are straightforward. If
$\psi=\Box\psi_0$, we have $X^\boxdot\vdash_\lgc K\psi_0\eq\psi_0^\boxdot$
by the induction hypothesis, hence
\[X^\boxdot\vdash_\lgc K\Box\psi_0\eq\Box\psi_0^\boxdot\]
using necessitation and~$\lgc K$. However,
\[X^\boxdot\vdash_\lgc K\Box\psi_0^\boxdot\eq\boxdot\psi_0^\boxdot\]
by~\ref{item:4}, whence
$X^\boxdot\vdash_\lgc K\Box\psi_0\eq\boxdot\psi_0^\boxdot=(\Box\psi_0)^\boxdot$.
\end{Pf*}
By \th\ref{cl:1} and the choice of~$\fii$, we have
$X^\boxdot\vdash_L\fii^\boxdot$, hence
\[\vdash_L\Box^{\le n}\ET X^\boxdot\to\fii^\boxdot\]
for some~$n$ by \th\ref{fact:dt}. (In fact, one can take for~$n$ the
modal degree of~$\fii$, but we will not need this.) In other words,
$\vdash_L\chi^\boxdot$, where
\[\chi=\Box^n\ET X\to\fii.\]
It remains to show that $\nvdash_{L_0}\chi$. The argument below
actually gives the slightly stronger conclusion $X\nvdash_{L_0}\fii$.

Since $\nvdash_{L_0}\fii$, there is a frame $\frm W\in C$
such that $\frm W\nmodel\fii$ by \th\ref{ass:dbl}. We can fix a model $\frm
M=\p{W,R,{\model}}$ based on~$\frm W$, and $w_0\in W$ such that
$\frm M,w_0\nmodel\fii$. We define a model~$2\frm
M=\p{2W,2R,{\model}}$ based on~$2\frm W$ by putting
\begin{align*}
2\frm M,\p{w,a}&\model q\iff a=1,\\
2\frm M,\p{w,a}&\model p_i\iff\frm M,w\model p_i
\end{align*}
for every $w\in W$, $a\in\{0,1\}$, and every variable~$p_i$ distinct
from~$q$. We have
\begin{equation}\label{eq:2}
2\frm M,\p{w,a}\model\psi\iff\frm M,w\model\psi
\end{equation}
for every $w\in W$, $a\in\{0,1\}$, and~$\psi\in\Sub(\fii)$ by
\th\ref{fact:modpres,obs:2w}, in particular $2\frm M,\p{w_0,a}\nmodel\fii$.
On the other hand, consider any formula
\[\Box(q^e\to\psi)\to\psi\]
from~$X$. If $2\frm M,\p{w,a}\nmodel\psi$, we have $\frm M,w\nmodel\psi$
by~\eqref{eq:2}, hence $\frm M,v\nmodel\psi$ for some~$w\mr Rv$ as $\frm
W\model L_0\Sset\lgc T$. (If $\frm W$ is a Kripke frame, it has to be
reflexive, in which case we can simply take $v=w$.) Then
$2\frm M,\p{v,e}\model q^e\land\neg\psi$ by~\eqref{eq:2}, and
$\p{w,a}\mr{2R}\p{v,e}$, thus $2\frm M,\p{w,a}\nmodel\Box(q^e\to\psi)$.

We have verified that $2\frm M\model X$ and~$2\frm M\nmodel\fii$. By
assumption, $2\frm W$ is an $L_0$-frame, hence $X\nvdash_{L_0}\fii$.
\end{Pf}

\begin{Cor}\th\label{cor:bdc}
The logics $\lgc T$, $\lgc{KTB}$, $\lgc{S4}$, $\lgc{S5}$, $\lgc{S4.2}$
and~$\lgc{S4.3}$ have the \sbc.
In particular, the boxdot conjecture is true.
\noproof\end{Cor}

It is natural to ask how far is \th\ref{ass:dbl} from being a
\emph{necessary} and sufficient criterion for (S)\wbc. On the positive
side, it is possible to show using similar techniques as in
Section~\ref{sec:motivation} that our condition gives full
characterization of (S)\wbc\ for well-behaved transitive logics:
\begin{Prop}\th\label{prop:char}
If $L_0\Sset\lgc{S4}$ has the finite model property, the following are
equivalent.
\begin{enumerate}
\item $L_0$ has the \wbc.
\item $L_0$ has the \sbc.
\item $L_0$ satisfies \th\ref{ass:dbl} for some~$C$.
\item $L_0$ satisfies \th\ref{ass:dbl} for every class~$C$ of general
frames with respect to which it is sound and complete.
\item $L_0$ is the smallest modal companion of some
superintuitionistic logic (see \cite[\S9.6]{cha-zax}).
\end{enumerate}
\end{Prop}
\begin{Pf}[sketch]
(ii)${}\to{}$(i) and (iv)${}\to{}$(iii) are trivial, and
(iii)${}\to{}$(ii) is \th\ref{thm:main}. (v)${}\to{}$(iv) follows from
the fact that $\frm W$ and~$2\frm W$ induce the same intuitionistic
frame.

(i)${}\to{}$(iii): Let $C$ be the class of all finite rooted $L_0$-frames.
For each $\frm F\in C$, let $\frm F^\I$ be the (not necessarily
transitive) frame obtained from~$\frm F$ by making all points
irreflexive, and let $L$ be the logic determined by $\{\frm F^\I:\frm
F\in C\}$. Since $(\frm F^\I)^\R=\frm F$, we have $\bdi L=L_0$, hence
$L\sset L_0$ by~(i). Since $L\Sset\lgc{wK4}=\lgc K\oplus\boxdot
p\to\Box\Box p$, we can use Fine frame formulas as in the transitive
case. In particular, $L\sset L_0$ implies that for every $\frm F\in
C$, there is a p-morphism $f\colon\frm G^\I\onto\frm F$ for some~$\frm
G\in C$. We use it to construct a mapping $g\colon G\to2F$ as
follows. If $x\in F$, and $c$ is a maximal cluster of~$\frm G$
intersecting~$f^{-1}[x]$, then $\lh{c\cap f^{-1}[x]}\ge2$ as $\frm
G^\I$ is irreflexive while $\frm F$ is reflexive. Thus, we can split
$c\cap f^{-1}[x]$ in two nonempty parts that are respectively mapped to
$\p{x,0}$ and~$\p{x,1}$ by~$g$. Other elements of~$f^{-1}[x]$ are
mapped, say, to~$\p{x,0}$. It is easy to check that $g$ is a
p-morphism of~$\frm G$ onto~$2\frm F$, hence $2\frm F$ is an
$L_0$-frame. Consequently, $L_0$ satisfies \th\ref{ass:dbl} wrt~$C$.

(iii)${}\to{}$(v): This can be shown using the machinery of
Zakharyaschev's canonical formulas (see \cite[\S9]{cha-zax} for an
explanation, which is outside the scope of this paper). Let
$\alpha(\frm F,D,\bot)$ be a canonical formula based
on a reflexive~$\frm F$, $\frm F'$ be a frame obtained from~$\frm F$ by
shrinking each cluster to (at least) half its original size, and $D'$
be the corresponding set of closed domains. Then any cofinal subreduction
of~$\frm W\model\lgc{S4}$ to~$\frm F'$ with the closed domain
condition for~$D'$ can be lifted to a cofinal subreduction of~$2\frm W$
to~$\frm F$ with CDC for~$D$. Using \th\ref{ass:dbl}, this shows that
$\vdash_{L_0}\alpha(\frm F,D,\bot)$ implies $\vdash_{L_0}\alpha(\frm
F',D',\bot)$. Repeating this process yields an axiomatization of~$L_0$ over~$\lgc{S4}$ by formulas~$\alpha(\frm F,D,\bot)$ where
$\frm F$ has only simple clusters, which means that $L_0$ is of the
form (v).
\end{Pf}
However, the situation appears to be more complicated in the case of
nontransitive~$L_0$.

\subsection*{Acknowledgement}
I would like to thank Alasdair Urquhart for introducing the problem to
me.

\bibliographystyle{mybib}
\bibliography{mybib}
\end{document}